\documentstyle[12pt]{amsart}

\let\Bbb\relax
\newfont{\Bb }{msbm10 scaled 1000} 
\newfont{\Bbb}{msbm10 scaled 1200}  
\newfam\eufam%
\font\euzw=eufm10 scaled 1200%
\font\euac=eufm9%
\textfont\eufam=\euzw\scriptfont\eufam=\euac%
\newcommand{\Z}{{\Bbb Z}}
\newcommand{\R}{{\Bbb R}}
\newcommand{\C}{{\Bbb C}}
\newtheorem{thm}{Theorem}[section] 
\newtheorem{prop}[thm]{Proposition} 
\newtheorem{cor}[thm]{Corollary} 
\newtheorem{lemma}[thm]{Lemma} 
\newtheorem{dof}{Definition}[section] 
\begin{document}
\title[On hyper K\"ahler manifolds associated to Lagrangean   
submanifolds]{On hyper K\"ahler manifolds associated to Lagrangean K\"ah\-ler 
submanifolds of $T^*\mbox{\C}^n$} 
\author{Vicente Cort\'es}
\thanks{Supported by the 
Alexander von Humboldt Foundation and MSRI (Berkeley).
Research at MSRI is supported in part by grant DMS-9022140.}

\address{\hskip-\parindent Mathematisches Institut der Universit\"at Bonn,
Beringstr. 1, 53115 Bonn, Germany}
\email{vicente@@msri.org {\rm or} vicente@@rhein.iam.uni-bonn.de} 

\curraddr{Mathematical Sciences
Research Institute\\
1000 Centennial Drive\\
Berkeley, CA 94720-5070}

\maketitle 
\pagenumbering{roman}
\begin{abstract}
For any Lagrangean K\"ahler submanifold $M \subset T^*\mbox{\Bb C}^n$, there
exists a canonical hyper K\"ahler metric on $T^*M$. A K\"ahler potential
for this metric is given by the generalized Calabi Ansatz of the 
theoretical physicists Cecotti, Ferrara and Girardello. 
This correspondence provides a method for the construction of 
(pseudo) hyper K\"ahler manifolds with large automorphism group. 
Using it, a class of pseudo hyper K\"ahler manifolds of complex
signature (2,2n) is constructed. For any hyper K\"ahler manifold $N$ in 
this class a group of automorphisms with a codimension one orbit on $N$ 
is specified. Finally, it is shown that the bundle of intermediate
Jacobians over the moduli space of gauged Calabi Yau 3-folds admits a 
natural pseudo hyper K\"ahler metric of complex signature (2,2n).
\end{abstract} 
\section*{Introduction}
The generalized Calabi Ansatz of Cecotti, Ferrara and Girardello
was discovered in the context of super string theory \cite{C-F-G}. 
Nevertheless, it provides a simple method for the 
construction of such classical geometric structures as hyper 
K\"ahler metrics. 

In the first part of this paper a self contained presentation of 
this construction is given. We explain how to a Lagrangean 
(pseudo) K\"ahler submanifold $M \subset T^*\mbox{\C}^n$ one canonically
associates a (pseudo) hyper K\"ahler metric on the complex symplectic 
manifold $T^*M$, s.\ Thm.~\ref{basicThm}. 

Then we study natural group actions on $M$ and on its 
cotangent bundle $T^*M$ preserving the special geometric structures. 
This opens the way for a systematic construction of (pseudo) hyper
K\"ahler manifolds of small cohomogeneity, s.\ Prop.~\ref{groupProp} and 
Cor.~\ref{codimCor}. 
Using the classification of certain
Lagrangean cones given in \cite{dW-VP} and \cite{C2}, we construct
examples of pseudo hyper K\"ahler manifolds of complex signature
(2,2n) admitting a group of automorphisms of cohomogeneity one, 
s.\ Thm.~\ref{examplesThm}. If $M$ is a cone of appropriate signature, 
then $M$ can be interpreted as formal moduli space of gauged 
Calabi Yau 3-folds, s.\ Prop.~\ref{formalProp} and Thm.~\ref{moduliThm}. 
In particular, the Lagrangean cones classified in  \cite{dW-VP} and \cite{C2} 
are models for moduli spaces with large automorphism group, s.\ 
Prop.~\ref{exformalProp} and Rem.\ 5. 

Finally, we prove that the pseudo  hyper K\"ahler structure on  the 
cotangent bundle $T^*M$ of a Lagrangean pseudo K\"ahler submanifold 
$M\subset T^*\mbox{\C}^n$ is always defined on a discrete fibre preserving 
quotient of $T^*M$ which is a torus bundle over $M$, s.\ 
Thm.~\ref{torusThm}. As a consequence, 
we obtain a natural pseudo hyper K\"ahler structure of complex signature 
(2,2n) on the bundle of intermediate Jacobians over the moduli space of 
gauged Calabi Yau 3-folds, where $n = h^{2,1}$, s.\ Thm.~\ref{JacThm}.  
 
\medskip\noindent 
{\bf Acknowledgements}\\
The author is very grateful to the Mathematical Sciences Research Institute,
especially to S.-S.\ Chern and R.\ Osserman, for hospitality and support.
He would also like to thank D.V.\ Alekseevsky, R.L.\ Bryant, A.B.\ Givental,
P.A.\ Griffiths, S.\ Katz, M.\ Kontsevich and A.\ Van Proeyen for 
discussions related to the subject of this paper.   
\newpage 
\tableofcontents  
\section[The canonical pseudo hyper K\"ahler metric on the cotangent bundle
$T^*M$ of a Lagrangean pseudo  K\"ahler submanifold 
$M\subset T^*\mbox{\C}^n$]{The canonical pseudo hyper K\"ahler metric on the 
cotangent bundle of a Lagrangean pseudo  K\"ahler submanifold 
$M\subset T^*\mbox{\C}^n$} 
\pagenumbering{arabic} 
Consider the following \label{fadP} {\bf fundamental algebraic data}: 
\begin{itemize} 
\item[1)] A complex symplectic vector space $(V, \omega )$, $\dim_{\mbox{\C}}V
= 2n$.  
\item[2)] A compatible real structure $\tau : V\rightarrow V$, i.e.\ a 
{\C}-antilinear involution such that the restriction of $\omega$ to its fix 
point set $V^{\tau}$ is a real symplectic structure.  
\end{itemize} 
Up to isomorphism, we can assume that $V = T^{\ast} \mbox{\Bbb C}^n$, 
$V^{\tau} = T^{\ast} \mbox{\Bbb R}^n$ and that $\omega$
is the standard symplectic structure of $T^{\ast} \mbox{\Bbb C}^n$, 
which is a real symplectic structure when restricted to 
$T^{\ast} \mbox{\Bbb R}^n$.  

Recall that a linear subspace of a (real or complex) symplectic 
vector space is called {\bf Lagrangean} if it is 
maximally isotropic and that a submanifold $M$ of a symplectic 
vector space is called {\bf Lagrangean
submanifold} if $T_mM$ is a Lagrangean subspace for every point $m\in M$. 

Given $(V, \omega , \tau)$ we can define a Hermitian form $\gamma$ of 
signature (n,n) on $V$ by 
\begin{equation} \gamma (u,v) = \sqrt{-1} \omega (u,\tau v)\, , \quad u,v\in V\, .\label{gammaEqu} \end{equation}
$(V,\gamma )$ is a pseudo K\"ahler manifold (of complex signature (n,n))
and hence the notion of pseudo K\"ahler submanifold is defined. 
In fact, a complex submanifold $M\subset V$ is called a {\bf pseudo K\"ahler
submanifold} of $(V,\gamma)$  if $\gamma |M$ is nondegenerate. 

Recall that a {\bf complex symplectic structure} on a complex manifold
is a holomorphic, closed and nondegenerate 2-form. 
\begin{dof} A (pseudo) K\"ahler manifold will be called
{\bf (pseudo) hyper K\"ahler manifold} if it admits a parallel
complex symplectic structure. \label{hyperDef} 
\end{dof} 

It follows from this definition that a K\"ahler manifold of complex 
dimension $2n$ is hyper K\"ahler if and only if the holonomy group 
(of the canonical connection) is a 
subgroup of $Sp(n) = U(2n)\cap Sp(n,\mbox{\C})$, cf.\ \cite{Be}. 
Here $Sp(n,\mbox{\C})$ denotes the symplectic group of $\mbox{\C}^{2n}$. 

\noindent 
{\bf Remark 1:} The complex symplectic structure $\omega$ defines on 
the (flat) pseudo K\"ahler manifold $(V, \gamma )$ the structure of pseudo 
hyper K\"ahler manifold. 
 
\begin{prop} \label{nondegProp} Given $(V, \omega , \tau)$ and $\gamma$ 
as above, a Lagrangean 
submanifold $M\subset (V, \omega )$ is a pseudo K\"ahler submanifold 
$M\subset (V,\gamma)$
if and only if $T_mM \cap \tau T_mM = 0$ for all $m\in M$. In particular,
$V^{\tau}\cap T_mM = 0$ is necessary.    
\end{prop} 

\noindent
{\bf Proof:} Let $L\subset  (V, \omega )$ be a Lagrangean subspace. Then 
$\gamma |L$ is nondegenerate if and only if $L\cap \tau L = 0$. $\Box$   

A connected Lagrangean pseudo  K\"ahler submanifold $M \subset (V, \omega , 
\tau)$ has a well defined complex signature $(k,l)$, $k+l = n$, namely the 
signature of the Hermitian form $\gamma |T_mM$, $m\in M$ arbitrary. 
If $l=0$, then $\gamma |T_mM$ is positively defined and $M$ is a K\"ahler
submanifold of $(V, \gamma )$. 

\begin{prop} \label{potMProp}For any pseudo  K\"ahler submanifold $M \subset (V, \gamma)$ 
the function $K^M(u):= \gamma (u,u)$, $u\in M$, is a pseudo K\"ahler 
potential. Any subgroup of $Aut (V,\gamma )$ which preserves $M$ 
acts on $M$ by holomorphic isometries. 
\end{prop}   

\noindent
{\bf Proof:} Let $u: U\rightarrow M$, $U\subset \mbox{\C}^m$, $m = 
\dim_{\mbox{\C}}M$, be a local holomorphic parametrization of $M$,
then 
\[ \frac{\partial^2 (K^M \circ u)}{\partial z^i\partial \overline{z}^j} \quad = \quad
\gamma ( \frac{\partial u}{\partial z^i},\frac{\partial u}{\partial z^j})\, .
\]
This proves the first claim. The second claim follows from the fact that
$Aut (V,\gamma ) \cong U(n,n)$ acts holomorphically and isometrically on the 
pseudo K\"ahler manifold $(V,\gamma )$ and that $M$ is a pseudo  K\"ahler 
submanifold of  $(V,\gamma )$. $\Box$ 

The linear automorphism group $Aut(V,\omega,\tau) \cong Aut(V^{\tau}, \omega |
V^{\tau}) \cong Sp(n,\mbox{\R})$ of our fundamental algebraic data acts
on the hyper K\"ahler manifold $(V,\gamma , \omega )$ by automorphisms,
i.e.\ by holomorphic isometries preserving the complex symplectic structure.
Here $Sp(n,\mbox{\R})$ denotes the real symplectic group in 2n variables.
 
\begin{dof} (cf.\ \cite{C-F-G}) The group 
\[ Aut_d (M) = \{ \varphi \in Aut(V,\omega,\tau)| \varphi M = M\} \]
is called the {\bf duality group} of the Lagrangean pseudo K\"ahler 
submanifold 
$M \subset (V, \omega , \gamma)$. 
\end{dof} 
The next proposition follows from Prop.~\ref{potMProp}.   
\begin{prop} \label{dualProp} The duality group $Aut_d(M)$ acts on $M$ 
holomorphically and isometrically. 
\end{prop}  

To define the generalized Calabi Ansatz of \cite{C-F-G}, we have to choose a 
{\bf Lagrangean splitting} for $V^{\tau}$, i.e.\ a decomposition 
\begin{equation} V^{\tau} = L_0 \oplus L_0' \label{LagsplEquI}\end{equation}
into two Lagrangean subspaces. A Lagrangean splitting for $V^{\tau}$ induces
a Lagrangean splitting for $V$: 
\begin{equation}  V = L \oplus L'\, , \label{LagsplEquII}\end{equation}
where $L$ and $L'$ are determined by 
\begin{equation} L=\tau L\, , \quad L'=\tau L' \quad \mbox{and} \quad L_0 = L^{\tau}\, ,
\quad L_0' = (L')^{\tau}\, .\label{LagsplEquIII}\end{equation}  
Moreover, we have canonical isomorphisms associated to (\ref{LagsplEquI})
and (\ref{LagsplEquII}):
\[ V^{\tau} \cong T^*L_0 \, , \quad V\cong T^*L\, .\] 

\begin{dof} \label{gpDef} A Lagrangean submanifold $M\subset (V,\omega)$ is 
in {\bf general position} with respect to the splitting  
$V = L \oplus L'$ if the 
projection $p: V\rightarrow L$ induces an isomorphism of $M$ onto its image.  
\end{dof}    

\begin{prop} \label{gpProp} A Lagrangean K\"ahler submanifold $M\subset 
(V,\omega ,\gamma )$
is in general position with respect to any splitting $V = L \oplus L'$ 
induced by a Lagrangean splitting of $V^{\tau}$, s.\ (\ref{LagsplEquI}),
(\ref{LagsplEquII}) and (\ref{LagsplEquIII}). 
\end{prop} 

\noindent
{\bf Proof:} It is sufficient to show that $T_mM\cap L' = 0$ for all 
$m\in M$. This follows from the fact that $\gamma$ is positively 
defined on $T_mM$ and zero on $L' = \tau L'$. $\Box$ 

Given fundamental algebraic data $(V,\omega , \tau )$,  
compatible Lagrangean splittings  
(\ref{LagsplEquI})--(\ref{LagsplEquIII}) and  a Lagrangean pseudo K\"ahler 
submanifold $M\subset  (V,\omega , \gamma )$ in general position with 
respect to
the given splitting, there is a {\bf canonical real structure}
$\rho '$ on $TM$. 
In fact, the projection $p: V\rightarrow L$ induces an isomorphism 
$T_mM\stackrel{\sim}{\rightarrow} L$ and hence we can define 
a real form of $(T_mM)^{\rho '}$ of $T_mM$ by the equation 
\[ dp\,(T_mM)^{\rho '} = L_0.\] 
Now we can define $\rho '$ (at the point m) as the {\C}-antilinear involution 
of $T_mM$ with fix point set $(T_mM)^{\rho '}$. 

We denote by $\rho$ the real structure on $T^*M$ which is dual to 
$\rho '$ and by $g^{-1}_m$ the Hermitian metric on $T_m^*M$ which
is inverse to $g_m = \gamma |T_mM$. Then the {\bf generalized Calabi 
Ansatz} of \cite{C-F-G} is given by the following potential on 
$T^*M$: 
\begin{equation} K(\sigma ) = K^M(\pi (\sigma )) + 
g^{-1}_{\pi (\sigma )} (\sigma + \rho (\sigma ),\sigma + \rho (\sigma ))\, ,
\quad \sigma \in T^*M\, , \label{hKpotEqu}\end{equation}
where $\pi : T^*M\rightarrow M$ is the natural projection and $K^M$ is 
the pseudo K\"ahler potential of $M$, s.\ Prop.~\ref{potMProp}. 

\begin{thm}\label{basicThm} Let $(V,\omega , \tau )$ be 
fundamental algebraic data, s.\ 
p.\ \pageref{fadP}, $\gamma$ the Hermitian form defined in (\ref{gammaEqu}),
$V = L \oplus L'$ a Lagrangean splitting as in 
(\ref{LagsplEquI})--(\ref{LagsplEquIII}) and $M\subset (V,\omega ,\gamma )$
a Lagrangean K\"ahler submanifold (resp.\ pseudo K\"ahler 
submanifold of complex signature $(k,l)$ in general position, s.\ Def.\ 
\ref{gpDef}).  Then the function $K$ on $T^*M$ associated
to these data, s.\ (\ref{hKpotEqu}), is the K\"ahler potential 
of a hyper K\"ahler metric $G$ on $T^*M$ (resp.\ the pseudo  K\"ahler 
potential of a pseudo hyper K\"ahler metric $G$ of complex 
signature $(2k,2l)$ on $T^*M$). More precisely, the standard complex 
symplectic 
structure $\Omega$ on the cotangent bundle $T^*M$ is parallel with respect
to the canonical connection of the K\"ahler (resp.\ pseudo K\"ahler) manifold
$N = (T^*M,G)$.      
\end{thm} 

\noindent
{\bf Remark 2:} In \cite{C-F-G} the correspondence $M \mapsto N$ is called 
{\it the c-map in rigid super symmetry}. There is also a {\it c-map in 
local super symmetry} $M \mapsto N$ as was proven in \cite{F-S}. In the latter
case, $N$ is a quaternionic K\"ahler manifold of negative Ricci curvature.  

\noindent
{\bf Proof:} First we derive the local coordinate expression for the 
field $G$ of Hermitian forms defined by the potential $K$ on $T^*M$. 
Using this expression, we show that $G$ is nondegenerate and has 
complex signature (2k,2l). Then we prove that the standard complex
symplectic structure $\Omega$ on $T^*M$ is parallel, also by a direct 
computation. 

Let us choose a linear isomorphism $L_0\cong \mbox{\R}^n$. 
Then we can identify $V^{\tau} = T^*\mbox{\R}^n$ and $V = T^*\mbox{\C}^n$. 
We denote by  $(q^1,\ldots ,q^n, p_1,\ldots ,p_n)$ the complex 
coordinates on $T^*\mbox{\C}^n$ which correspond to the 
standard coordinates on $\mbox{\R}^n$. In these coordinates 
$\omega = \sum_{i=1}^n dq^i\wedge dp_i$. 

Since $M$ is a Lagrangean submanifold in general position,
s.\ Prop.\ \ref{gpProp},  it is the image of a closed and hence locally 
exact section of $T^*\mbox{\C}^n$. In  other words, we can describe $M$ 
locally by equations of the form
\begin{equation} p_i = \frac{\partial F (q^1,\ldots ,q^n)}{\partial q^i}\, , \quad 
i = 1,\ldots , n \, , \label{LagEqu} \end{equation} 
where $F(q^1,\ldots q^n)$ is a locally defined holomorphic function 
of n variables. Remark that for any locally defined holomorphic function F
the equations (\ref{LagEqu}) define a Lagrangean submanifold in general 
position, namely the image of the exact section $dF$ of $T^*\mbox{\C}^n$.  
 
Denote by $z^i:= q^i|M$, $i = 1,\ldots ,n$, the natural complex coordinates
on $M$ and by $(z^1,\ldots ,z^n, w_1,\ldots ,w_n)$ the corresponding  
complex coordinate system for $T^*M$. To unify the notation put 
$z^{i'} := w_i$ and $(z^I) := ((z^i),(z^{i'}))$. Recall that the (pseudo)
K\"ahler metric on $M$ is $g = \gamma |M$, so in our coordinates $(z^i)$:
 \begin{equation} \label{gijEqu} g_{ij} = g_{ji} = 
\gamma (\frac{\partial}{\partial z^i},
\frac{\partial}{\partial z^j}) = \sqrt{-1}(\overline{F}_{ij}- F_{ij})\, ,
\end{equation} 
where $F_i = p_i|M$, $F_{ij} = \frac{\partial F_i}{\partial z^j}$, 
$F_{ijk} = \frac{\partial F_{ij}}{\partial z^k}$ etc.\ and we have
used that 
\[ \frac{\partial}{\partial z^i} \quad = \quad  \frac{\partial}{\partial q^i}
+ \sum_j F_{ij} \frac{\partial}{\partial p_j}\, .\] 
The Hermitian fibre metric $g^{-1} = (g^{ij})$ on $T^*M$ is defined by the 
equation
\begin{equation}\label{g^ijEqu} \sum_j g^{ij}g_{jk} = \delta_k^i\, .
\end{equation} 
With respect to the coordinates $(z^1,\ldots , z^n,w_1,\ldots ,w_n)$ on 
$T^*M$ the potential $K$ reads 
\begin{equation} K(z^1,\ldots , z^n,w_1,\ldots ,w_n) = 
K^M(z^1,\ldots , z^n) + 
\sum_{ij} g^{ij} (w_i + \overline{w}_i)(w_j + \overline{w}_j)\, ,
\end{equation}  
where 
\begin{equation} \label{potcoordEqu} K^M(z^1,\ldots , z^n) 
= \gamma ({dF|}_{(z^1,\ldots , z^n)},
{dF|}_{(z^1,\ldots , z^n)}) = \sqrt{-1} \sum_{ij} (z^i\overline{F}_i -   
\overline{z}^iF_i)\, . 
\end{equation}  
Now we compute $G = (G_{IJ})$, where 
\[ G_{IJ} = \frac{\partial^2 K}{\partial z^I \partial \overline{z}^J}\, ,\quad
I,J \in\{ 1,\ldots ,n, 1',\ldots , n'\} \, .\] 
In the following, lower indices separated
by komma will always denote partial derivatives, e.g.\ $g_{ij,k\bar{l}m'} =
\frac{\partial^3g_{ij}}{\partial z^k \partial \overline{z}^l \partial w_m}$. 
The equations 
\begin{eqnarray*} G_{ij}  & =  & g_{ij} + \sum_{kl}{g^{kl}}_{,i\bar{j}}
(w_k + \overline{w}_k)(w_l + \overline{w}_l)\\
& = & g_{ij} + 2\sum_{krspql}g^{kr}g^{sp}g^{ql}F_{rsi}\overline{F}_{pqj}
(w_k + \overline{w}_k)(w_l + \overline{w}_l)\\
G_{ij'} & = & \overline{G_{j'i}} = 2\sqrt{-1}\sum_{kpq}g^{kp}F_{pqi}g^{qj}
(w_k + \overline{w}_k)\\
G_{i'j'} & = & 2g^{ij} 
\end{eqnarray*} 
follow immediately from the basic formulas
\begin{eqnarray}
{g^{kl}}_{,i} & = & - \sum_{pq}g^{kp} g_{pq,i}g^{ql} 
= \sqrt{-1}\sum_{pq}g^{kp} F_{pqi}g^{ql}\label{g-1'EquI}\\
{g^{kl}}_{,\bar{j}} & = & - \sum_{pq}g^{kp} g_{pq,\bar{j}} g^{ql} 
= -\sqrt{-1}\sum_{pq}g^{kp} \overline{F}_{pqj}g^{ql}\, .\label{g-1'EquII}
\end{eqnarray} 
If we put $b = (b^i_j)_{i,j = 1,\ldots , n}$, $b^i_j := G_{ji'}$, then 
\[ G = (G_{IJ}) \quad = \quad \left( \begin{array}{ll} 
g + \frac{1}{2}b^tg\overline{b} & b^t \\
\overline{b} & 2 g^{-1} \end{array} \right) \]
and we can easily invert this matrix: 
\[ G^{-1} = (G^{IJ}) \quad = \quad \left( \begin{array}{ll} 
g^{-1} & -\frac{1}{2} b\\
-\frac{1}{2} \overline{b}^t & \frac{1}{2}(g + \frac{1}{2}\overline{b}^tgb)
\end{array} \right) \, .\] 
This shows that $G$ is a pseudo K\"ahler metric on $T^*M$. In particular,
$(T^*M,G)$ has a well defined signature over each connected component 
of $M$. We may assume that $M$ is connected and $(M,g)$ has complex signature
(k,l). Then $(T^*M,G)$ has signature (2k,2l) near the zero section 
$M\subset T^*M$ (i.e.\ for $b\rightarrow 0$) and hence everywhere. 

Now we show that $N = (T^*M,G)$ is a pseudo 
hyper K\"ahler manifold, s.\ Def.~\ref{hyperDef}. The complex manifold
$T^*M$ has the canonical complex symplectic structure 
$\Omega = \sum dz^i\wedge dw_i$. We will show that $\nabla \Omega = 0$ for 
the covariant derivative $\nabla$ of the pseudo K\"ahler manifold $N$. 
Let us denote by $\Gamma^I_{JK}$ the Christoffel 
symbols of the pseudo K\"ahler metric $G$ on $N$. 
We recall that $\Gamma^I_{JK} = \sum_LG^{LI}G_{JL,K}$, s.\ e.g.\ 
\cite{K-N2} for the basic theory of K\"ahler manifolds.  
        
\begin{lemma} \label{ChristoffelLemma} The complex symplectic structure $\Omega$ of $T^*M$ is parallel
with respect to the pseudo K\"ahler metric $G$ if and only if the Christoffel
symbols $\Gamma^I_{JK}$ have the following symmetries.
\begin{itemize}
\item[(i)] $\Gamma^i_{Jk} = - \Gamma^{k'}_{Ji'}$,
\item[(ii)] $\Gamma^i_{Jk'} = \Gamma^{k}_{Ji'}$, $\Gamma^{i'}_{Jk} = 
\Gamma^{k'}_{Ji}$. 
\end{itemize} 
\end{lemma}

\noindent
{\bf Proof:} It is straightforward to check that equations (i) and (ii)
are equivalent to 
\[ \nabla_{\frac{\partial}{\partial z^l}}  \Omega = 0\, , \quad 
l = 1,\ldots , n\, ,\] 
if $J = j = 1,\ldots , n$ and equivalent to 
\[ \nabla_{\frac{\partial}{\partial w_l}}  \Omega = 0\, , \quad 
l = 1,\ldots , n\, ,\]
if $J = j' = 1',\ldots , n'$. $\Box$ 

The equations (i) and (ii) of Lemma \ref{ChristoffelLemma} are verified 
by a direct computation of the Christoffel symbols 
$\Gamma^I_{JK}$. 
This finishes the proof of Thm.~\ref{basicThm}. $\Box$ 

\noindent 
{\bf Remark 3:} Instead of Lemma \ref{ChristoffelLemma} one can also
use \cite{H} Lemma 6.8, cf.\ \cite{C-F-G}. We remark that to 
the pseudo hyper K\"ahler manifold $(T^*M,G,\Omega )$ constructed 
in Thm.\ \ref{basicThm} we can canonically associate a parallel 
hypercomplex structure $(J_1,J_2,J_3)$, which is Hermitian with 
respect to the pseudo Riemannian metric 
$\langle \cdot ,\cdot \rangle = {Re}\, G$. Here $J_1$ is the standard complex
structure of the (holomorphic) cotangent bundle $T^*M$; the complex 
structures $J_2$ and $J_3$ are defined by the equation 
\[ \Omega (v,w) \; = \; \langle J_2 v,w\rangle + 
\sqrt{-1}\, \langle J_3 v,w\rangle \, , \quad v,w \in T(T^*M)\, .\]    

\noindent  
{\bf Example 1:} The simplest example of Lagrangean pseudo K\"ahler
submanifold $M\subset T^*\mbox{\C}^n$ is a Lagrangean subspace $L$ such that
$L\cap \tau L = 0$. If e.g.\ 
\[ L =  {span}_{\mbox{\C}}\{ 
\frac{\partial}{\partial q^1} + i\frac{\partial}{\partial p_1},
\ldots , 
\frac{\partial}{\partial q^k} + i\frac{\partial}{\partial p_k},
\frac{\partial}{\partial q^{k+1}} 
- i\frac{\partial}{\partial p_{k+1}}, \ldots ,
\frac{\partial}{\partial q^n} 
- i\frac{\partial}{\partial p_n}\}\, ,\]
$i= \sqrt{-1}$, then $\gamma |L$ has complex signature (k,n-k) and $(T^*L,G)$ is the 
flat model of pseudo hyper K\"ahler manifold of quaternionic signature
(k,n-k), i.e.\ of complex signature (2k,2n-2k).    
\section{A class of pseudo hyper K\"ahler manifolds admitting a group
of au\-to\-mor\-phisms of cohomogeneity one} \label{exSec} 
Let $(V,\omega , \tau )$ be our fundamental algebraic data, s.\ 
p.\ \pageref{fadP}, $\gamma$ the Hermitian form defined in (\ref{gammaEqu}),
$V = L \oplus L'$ a Lagrangean splitting as in 
(\ref{LagsplEquI})--(\ref{LagsplEquIII}) and $M\subset (V,\omega ,\gamma )$
a Lagrangean pseudo K\"ahler 
submanifold in general position, s.\ Def.\ 
\ref{gpDef} and Prop.\ \ref{gpProp}. By Thm.\ \ref{basicThm} to these data we
canonically associate the pseudo hyper K\"ahler manifold $N = (T^*M,G)$
with standard complex symplectic structure $\Omega$. 

Now we are  interested in natural group actions on $M$ and $N$ preserving
the given geometric structures. Recall that the duality group $Aut_d(M)$ 
of $M$ consists of those linear automorphisms of  
$(V,\omega , \tau )$ which preserve $M$. It acts holomorphically 
and isometrically on the pseudo K\"ahler manifold $(M,g)$, s.\ 
Prop.~\ref{dualProp}. Let us consider the subgroup $Aut_{sd}(M)$ 
of $Aut_d(M)$ preserving 
also the Lagrangean splittings (\ref{LagsplEquI})--(\ref{LagsplEquIII}), i.e.\ 
\[ Aut_{sd}(M) \;  = \; \{ \varphi \in Aut_d(M)| \varphi L = L\, , 
\; \varphi L' = L'\}\, .\] 
Since $Aut_{sd}(M) \subset Aut_d(M) \hookrightarrow 
Sp(V^{\tau}, \omega |V^{\tau})\subset 
GL(V^{\tau}) \subset Aff(V^{\tau}) = V^{\tau} \mbox{\Bbb o} GL(V^{\tau})$,
we can also consider the  corresponding affine group
\[ V^{\tau} \mbox{\Bbb o} Aut_{sd}(M) \subset V^{\tau} \mbox{\Bbb o} Aut_d(M)
\hookrightarrow Aff(V^{\tau})\, .\] 


We will show that the affine group $V^{\tau} \mbox{\Bbb o} Aut_{sd}(M)$ acts 
on $N$ preserving the 
pseudo hyper K\"ahler structure. First we define a fibre preserving 
action of the vector group $V^{\tau} \cong \mbox{\R}^{2n}$ on 
$T^*M$. For this we consider the quotient map $V \rightarrow V/T_mM$. 
Using the inclusion $V^{\tau}\subset V$ and the isomorphism 
$V/T_mM \cong T^*_mM$ given by the symplectic structure $\omega$, 
the quotient map induces an isomorphism of real vector spaces, cf.\ 
Prop.~\ref{nondegProp}:  
\begin{equation}\label{psimEqu} \psi_m :V^{\tau} \rightarrow T^*_mM \, .
\end{equation}  
Now we define the action of an element $v\in V^{\tau}$ on 
$T^*M$ by 
\begin{equation} \label{vectoractionEqu} 
T^*_mM \ni \sigma \; \mapsto \; \sigma + \sqrt{-1} \psi_m(v) \in T^*_mM\, .
\end{equation} 
Choosing linear coordinates $(v^1,\ldots ,v^n)$ for $L_0\cong \mbox{\R}^n$
induces coordinates $(v_1, \ldots ,v_n)$ for $L_0' \cong L^*_0$ 
and $(z^1,\ldots ,z^n, w_1,\ldots ,w_n)$ for $T^*M$ as in the 
previous section. In these coordinates the action of 
$v = (v^1,\ldots ,v^n,v_1, \ldots ,v_n)$ on $T^*M$ defined in 
(\ref{vectoractionEqu}) reads (cf.\ \cite{C-F-G}):
\[ (z^1,\ldots ,z^n, w_1,\ldots ,w_n) 
\quad \mapsto \quad 
 (\tilde{z}^1,\ldots ,\tilde{z}^n, \tilde{w}_1,\ldots ,\tilde{w}_n)\, 
\]    
\begin{equation} \label{coordEqu}\tilde{z}^i \quad = \quad z^i\, , \qquad 
\tilde{w}_i \quad = \quad 
w_i - \sqrt{-1}(v_i - \sum_j F_{ij}v^j)\, ,\end{equation}  
where we recall that $F_{ij} = \frac{\partial^2 F}{\partial z^i \partial z^j}$
are the second derivatives of the local holomorphic function 
$F(z^1,\ldots ,z^n)$ defining the Lagrangean submanifold $M$, s.\ 
(\ref{LagEqu}). 

\begin{prop} \label{groupProp} The group $V^{\tau} \mbox{\Bbb o} Aut_{sd}(M)
\hookrightarrow Aff(V^{\tau})$ acts naturally on $N$ by automorphisms of the 
pseudo
hyper K\"ahler structure (s.\ Thm.~\ref{basicThm}), i.e.\ by holomorphic 
isometries preserving the 
complex symplectic structure $\Omega$. The action of $Aut_{sd}(M)$ is by point 
transformations of $N = T^*M$ and that of $V^{\tau}$, given by 
(\ref{vectoractionEqu}), is fibre preserving and simply transitive on each
fibre. 
\end{prop} 

\noindent
{\bf Proof:} $Aut_d(M)$ and thereby its subgroup $Aut_{sd}(M)$ acts 
holomorphically on $M$ and hence by complex symplectomorphisms on 
$T^*M$, namely by holomorphic point transformations. Moreover,    
$Aut_{sd}(M)$ preserves the pseudo K\"ahler potential (\ref{hKpotEqu})
on $T^*M$ associated to the choice of Lagrangean splitting (\ref{LagsplEquI}).
This shows that $Aut_{sd}(M)$  acts on $N$ by automorphisms of the 
pseudo hyper K\"ahler structure. 

Consider the map 
\[ \psi : M\times V^{\tau} \rightarrow T^*M\, , \quad \psi (m,v) := \psi_m(v)
\, .\] 
>From the fact that $\psi$ is $Aut_{d}(M)$-equi\-va\-riant it follows that
(\ref{vectoractionEqu}) extends the action of $Aut_d(M)$ on $N$ to a 
holomorphic action of $V^{\tau}\mbox{\Bbb o} Aut_{sd}(M)$ on $N$. 
We check that the action of $v\in V^{\tau}$, s.\ (\ref{coordEqu}), preserves 
the complex symplectic structure $\Omega = \sum_i dz^i \wedge dw_i$:
\[ \sum_i (d\tilde{z}^i \wedge d\tilde{w}_i -  dz^i \wedge dw_i) =
\sqrt{-1}\sum_{ij} v^j dz^i\wedge dF_{ij} = \sqrt{-1}\sum_{ijk} v^j F_{ijk} 
dz^i\wedge dz^k = 0\, .\]
Next we check that under the action of $v\in V^{\tau}$ the pseudo 
K\"ahler potential $K$ of $N$, s.\ (\ref{hKpotEqu}), changes only by
a pluriharmonic function. Using (\ref{gijEqu})-(\ref{potcoordEqu}) and 
(\ref{coordEqu}) we obtain: 
\[ K(\tilde{z}^1,\ldots ,\tilde{w}_n) 
- K(z^1,\ldots ,w_n) 
= -2\sum_j v^j (w_j + \overline{w}_j) + \sum_{ij} g_{ij}v^iv^j\, . \]
The result is pluriharmonic, since the $v^i$ are constants and $g_{ij}$
is pluriharmonic as sum of a holomorphic and an antiholomorphic function,
s.\ (\ref{gijEqu}). Now it only remains to show that $V^{\tau}$ acts 
simply transitively on each fibre of $T^*M \rightarrow M$. This is clear since
(\ref{psimEqu}) is an isomorphism of real vector spaces. $\Box$ 

\begin{cor} \label{codimCor} If a subgroup $A \subset Aut_{sd}(M)$ has an 
orbit of codimension r 
on $M$ then the subgroup $V^{\tau} \mbox{\Bbb o} A \subset V^{\tau} 
\mbox{\Bbb o} Aut_{sd}(M)$ has an orbit of codimension r on $N$.
\end{cor} 

Our aim is to use Cor~\ref{codimCor} for the construction of pseudo hyper 
K\"ahler manifolds $N$ with the smallest possible cohomogeneity of the 
group $V^{\tau} \mbox{\Bbb o} Aut_{sd}(M)$. 
(Recall that the cohomogeneity of a Lie group acting on 
a manifold is the minimal codimension of its orbits.) 

The pseudo K\"ahler
potential $K^M(u) = \gamma (u,u)$, $u\in M$, defines an   
$Aut_d(M)$-in\-va\-riant function on $M$ and hence $K^M \circ \pi$,
$\pi : T^*M \rightarrow M$ the projection, is an 
$V^{\tau} \mbox{\Bbb o} Aut_{d}(M)$-in\-va\-riant function on $N$. This function cannot be constant 
on an open set, since $K^M$ is the potential of a (nondegenerate) metric. 
Therefore, $V^{\tau} \mbox{\Bbb o} Aut_{d}(M) \supset V^{\tau} \mbox{\Bbb o} 
Aut_{sd}(M)$ has no open orbit on $N$ and hence is of co\-ho\-mo\-ge\-nei\-ty
at least one. 

Next we will present a class of pseudo hyper K\"ahler manifolds $N$ for 
which the cohomogeneity of $V^{\tau} \mbox{\Bbb o} 
Aut_{sd}(M) \subset V^{\tau} \mbox{\Bbb o} Aut_{d}(M)$ is in fact one. 
This class is associated, 
via the correspondence of Thm.~\ref{basicThm},  to an interesting class
of Lagrangean pseudo K\"ahler submanifolds $M\subset T^*\mbox{\C}^{n+1}$  
in general position (with respect to the standard Lagrangean splitting
$T^*\mbox{\C}^{n+1} = \mbox{\C}^{n+1} \oplus (\mbox{\C}^{n+1})^*$ into 
positions and momenta). The latter class can be defined by the following 
additional conditions\label{i-iiiDef}  
\begin{itemize}
\item[(i)]  \label{i-iii} $M$ is a cone, i.e.\ $\lambda M = M$ for all 
$\lambda \in \mbox{\C}-\{ 0\}$. 
\item[(ii)] The third fundamental form of $M\subset T^*\mbox{\C}^{n+1}$ is 
given by a homogeneous cubic polynomial $h(x^1,\ldots ,x^n)$ with real 
coefficients, s.\ remarks below.  
\item[(iii)] The real hypersurface $\{ h = 1\}\subset \mbox{\R}^n$ 
admits an open orbit ${\cal H} \subset \{ h = 1\}$ of a subgroup of 
$GL(n,\mbox{\R})$. Moreover, the second fundamental form of $\cal H$ is 
negatively defined. 
\end{itemize}   
The complete classification of such manifolds $M$ was first obtained in
\cite{dW-VP}; an alternative, conceptual approach was developped
in \cite{C2}.  

First of all let us explain the meaning of (ii).  The notion of 3rd 
fundamental form of a Lagrangean cone was introduced by Bryant and 
Griffiths \cite{B-G} in the study of period maps for general
Calabi Yau 3-folds, s.\ Section \ref{CYSec}. The 
{\bf third fundamental form} $\theta$ of a 
Lagrangean submanifold $M\subset T^*\mbox{\C}^{n+1}$ in general 
position satisfying (i) is by definition the 3rd fundamental form
of the immersion 
\[ \phi : (q^1,\ldots ,q^n) \; \mapsto \; (\pi |M)^{-1}(1,q^1,\ldots ,q^n)\, 
,\]
where $\pi : T^*\mbox{\C}^{n+1} \rightarrow \mbox{\C}^{n+1}$ is the 
natural projection. More precisely, $\theta$ is the section of 
$S^3T^*\mbox{\C}^n$  whose value $\theta_q 
\in S^3T^*_q\mbox{\C}^n\cong S^3(\mbox{\C}^n)^*$ at 
$q = (q^1,\ldots ,q^n)$, $(1,q)\in \pi (M)$, is given by the cubic form   
\[ \theta_q (x^1,\ldots ,x^n)\; = \; \sum_{ijk=1}^n x^ix^jx^k
\frac{\partial^3 \phi (q)}{\partial q^i \partial q^j \partial q^k}\, .\] 
So the condition (ii) is satisfied if and only if the homogeneous cubic 
polynomial $\theta_q$ is independent of $q$, i.e.\ 
$\theta_q(x^1,\ldots ,x^n) = h(x^1,\ldots ,x^n)$,  and has 
(constant) real coefficients 
$\frac{\partial^3 \phi (q)}{\partial q^i \partial q^j \partial q^k}\in 
\mbox{\R}$. 
   
In the following discussion we recall all needed basic facts 
about the class of Lagrangean submanifolds defined above; for 
complete details the reader is referred to \cite{dW-VP} and \cite{C2}.   
Under the condition (i), we can consider the complex manifold 
$P(M) \subset P(T^*\mbox{\C}^{n+1}) 
\cong P_{\mbox{\C}}^{2n+1}$, which is a projectivized cone. 
If the potential $K^M$ does not vanish on $M$, then $P(M)$ admits a canonical 
Hermitian form $g^{P(M)}$ invariant under the projective action
of the duality group $Aut_d(M)$, cf.\ (\ref{sKEqu}). 
(The corresponding potential is 
$\log |K^M|$.\label{sKP}) Under the conditions (ii) and (iii), there exists an 
open 
subcone ${\cal C} \subset M\subset T^*\mbox{\C}^{n+1}$ such that 
${Aut}_{sd}({\cal C})$ acts transitively on  $P({\cal C})\subset P(M)$,
the Hermitian form $\gamma$ has complex signature (1,n) on ${\cal C}$ and 
$-g^{P({\cal C})}$ 
is a (positively defined) K\"ahler metric on $P({\cal C})$. 
Without restriction of generality we assume $M = {\cal C}$.  From the 
preceding remarks it follows that  
${Aut}_{sd}(M)$ has an orbit of codimension one on $M$ and by 
Prop.~\ref{groupProp} and Cor.~\ref{codimCor} the group 
$V^{\tau} \mbox{\Bbb o} Aut_{sd}(M)$, $V^{\tau} = T^*\mbox{\R}^{n+1} 
\cong \mbox{\R}^{2n+2}$,  of automorphisms of the 
hyper K\"ahler manifold $N$ has cohomogeneity one. 

Now we describe the group ${Aut}_{sd}(M)$ in more detail. Remark that the 
subgroup
of $Sp(n+1,\mbox{\R}) \subset GL(2n+2,\mbox{\R})$ preserving the standard
Lagrangean splitting $T^*\mbox{\R}^{n+1} = \mbox{\R}^{n+1} \oplus 
(\mbox{\R}^{n+1})^*$ is the group of linear point transformations of 
$T^*\mbox{\R}^{n+1}$, which is canonically isomorphic to $GL(n+1,\mbox{\R})$. 
This gives rise to the embedding 
\[ \iota : Aut_{sd}(M) \hookrightarrow GL(n+1,\mbox{\R})\, , \quad 
\varphi \mapsto \varphi |\mbox{\R}^{n+1}\, ,\] 
where $\mbox{\R}^{n+1} \subset T^*\mbox{\R}^{n+1}$ is the zero section. 
$GL(n+1,\mbox{\R})$ contains the affine group $Aff(n,\mbox{\R}) = 
\mbox{\R}^n \mbox{\Bbb o} GL(n,\mbox{\R})$ as a subgroup: 
\[ Aff(n,\mbox{\R}) \ni (v,B) \; \mapsto \; 
\left( \begin{array}{cc} 1 & 0\\
v & B\end{array} \right) \in GL(n+1,\mbox{\R}) \, .\]   
So we can define the subgroup 
\[ A(M) \; := \; \iota^{-1} (Aff(n,\mbox{\R})) \subset Aut_{sd}(M)\, ,\]
which is embedded in the affine group $Aff(n,\mbox{\R})$ via the map $\iota$.

By construction, s.\ \cite{dW-VP} and \cite{C2}, the class of submanifolds
$M\subset T^*\mbox{\C}^{n+1}$ defined above has the following property: 
The group $A(M)$ has an open orbit on $P(M)$ and (by restricting to an 
open subcone of $M$ if necessary) we can assume that $A(M)$ acts transitively
on $P(M)$. Moreover, it was proven in \cite{C2} that $(P(M), -g^{P(M)})$ 
is isomorphic, as $A(M)$-K\"ahler manifold, to a K\"ahlerian Siegel 
domain $(U,g)$ of first type. The group $\iota (A(M)) \cong A(M)$ acts 
naturally on $(U,g)$ by affine transformations which are (holomorphic)
isometries for the K\"ahler metric $g$.  

The complete list of K\"ahlerian Siegel domains $(U,g)$ corresponding to 
the class of submanifolds $M\subset T^*\mbox{\C}^{n+1}$ defined above, 
s.\ (i)-(iii), was given in \cite{C2} Thm.~2.8. For convenience of 
the reader, we recall that $U$ has rank 2 or 3. The rank 2 domains in the list 
are numerated by the nonnegative integers and the rank 3 domains by special
isometric maps or, equivalently, by ${\mbox{\Z}}_2$-graded Clifford modules 
(up to equivalence defined in \cite{C1} Def.\ I.10 and Def.\ II.5, cf.\ 
Def.\ II.4, Prop.\ II.25 and Prop.~II.26). Recall that a 
${\mbox{\Z}}_2$-{\bf graded Clifford module of order k} is a 
${\mbox{\Z}}_2$-graded module $\Psi = \Psi_0 \oplus  \Psi_1$ of the 
real Clifford algebra ${C\! \ell}_k$, $k = 0,1,2,\ldots$, s.\ e.g.\ 
\cite{L-M}.      Summarizing our discussion and applying Cor.~\ref{codimCor}
to $A = A(M)$, we obtain the following theorem. 

\begin{thm} \label{examplesThm} To any $p\in \{ 0,1,2,\ldots\}$ (resp.\ 
${\mbox{\Z}}_2$-graded Clifford module $\Psi$ of order $k \in 
\{ 0,1,2,\ldots\}$) we can canonically associate a Lagrangean 
pseudo K\"ahler submanifold in general position $M(p) \subset 
T^*\mbox{\C}^{n+1}$, $n = 2+p$, (resp. $M(\Psi)\subset 
T^*\mbox{\C}^{n+1}$, $n = k+3+ \dim \Psi$) satisfying conditions (i)-(iii) 
on p.\ \pageref{i-iii}. Conversely, any Lagrangean pseudo K\"ahler 
submanifold of $T^*\mbox{\C}^{n+1}$ in general position satisfying 
(i)-(iii) contains one of the manifolds $M$ above (i.e.\ $M = M(p)$ or 
$M = M(\Psi )$) as an open subcone. 
For the manifolds $M$ above, the affine group $A(M) 
\stackrel{\iota}{\hookrightarrow} Aff(n,\mbox{\R})$ acts transitively
and isometrically on the K\"ahler manifold $(P(M), -g^{P(M)})$ 
by projective linear transformations. Moreover, $(P(M), -g^{P(M)})$
is isomorphic as $A(M)$-K\"ahler manifold  to a K\"ahlerian Siegel domain
of type I with transitive affine action of $\iota (A(M)) 
\subset Aff(n,\mbox{\R})$. 

The pseudo hyper K\"ahler manifold $N = (T^*M,G)$ associated to any of these
manifolds $M$ by Thm.~\ref{basicThm} has complex signature $(2,2n)$. 
Finally, the affine group $\mbox{\R}^{2n+2} \mbox{\Bbb o} A(M) 
\subset Aff(2n+2,\mbox{\R})$ acts on $N$ by automorphisms of the pseudo
hyper K\"ahler structure with an orbit of codimension one. 
\end{thm} 

\noindent 
{\bf Remark 4:} The manifolds $(T^*M,G)$, $M = M(p)$ or $M = M(\Psi )$,
should be thought of as natural pseudo hyper K\"ahlerian versions of 
Alekseevsky's homogeneous quaternionic K\"ahler manifolds, cf.\ \cite{A},
\cite{Ce}, \cite{dW-VP}, \cite{dW-V-VP}, \cite{C1}, \cite{A-C} and
\cite{C2}.             
\section{The pseudo hyper K\"ahler metric of the bundle of intermediate
Jacobians over the moduli space of gauged Calabi Yau 3-folds} 
\label{CYSec} 
Let $(V,\omega ,\tau )$ be our fundamental algebraic data, s.\ p.\ 
\pageref{fadP},
$\gamma$ the Hermitian form defined in (\ref{gammaEqu}), $V = L \oplus L'$ a
Lagrangean splitting as in (\ref{LagsplEquI})-(\ref{LagsplEquIII}), 
$M\subset (V,\omega ,\gamma )$ a Lagrangean pseudo K\"ahler submanifold
in general position, s.\ Prop.\ \ref{gpProp}, and, finally, $\Gamma \subset
V^{\tau}$ a (cocompact) lattice. Using the isomorphism $V/T_mM\cong T_m^*M$
induced by the symplectic form $\omega$, we can identify $T^*M$ with the 
normal bundle ${\cal N} \rightarrow M$ of $M$ in $V$, ${\cal N}_m = V/T_mM$. 
Since $V^{\tau} \supset \Gamma$ has zero intersection with $T_mM$, s.\ Prop.\ 
\ref{nondegProp}, $\Gamma$ projects to a lattice ${[}\Gamma{]}$ in 
${\cal N}_m = V/T_mM$ and ${\cal N}_m/ {[}\Gamma{]}$ is a complex 
torus. Let us denote by ${\cal N}/\Gamma$ the corresponding (holomorphic)
torus bundle over $M$.  Remark that via the isomorphism ${\cal N}_m 
\cong T_m^*M$ the lattice ${[}\Gamma{]}\subset {\cal N}_m$ corresponds to 
the lattice $\psi_m(\Gamma ) \subset T^*_mM$, s.\ (\ref{psimEqu}), and we can 
identify ${\cal N}/\Gamma$ with the quotient of $T^*M$ by the action of 
$\Gamma \subset V^{\tau}$ defined in equation (\ref{vectoractionEqu}). 

\begin{thm} \label{torusThm} The pseudo hyper K\"ahler structure on 
$T^*M \cong {\cal N}$ 
constructed in Thm.\ \ref{basicThm} induces a pseudo  hyper K\"ahler 
structure on the torus bundle $T^*M/ \Gamma \cong{\cal N}/\Gamma$ 
for any lattice $\Gamma \subset V^{\tau}$. 
\end{thm}

\noindent
{\bf Proof:} It was proven in Prop.\ \ref{groupProp} that the action 
of $V^{\tau} \supset \Gamma$ preserves the complex symplectic structure
$\Omega$ and the pseudo hyper K\"ahler metric $G$ on $T^*M$. $\Box$ 

The purpose of this section is to use Thm.\ \ref{torusThm} for the 
construction of a pseudo hyper K\"ahler 
structure on the bundle of intermediate Jacobians over the moduli
space of gauged Calabi Yau 3-folds. For this we have to review
some known facts about the moduli space, thereby relating it to 
Thm.\ \ref{torusThm}. 

Let $X$ be a {\bf (general) Calabi Yau 3-fold}, i.e.\ a compact K\"ahler 
3-fold with holonomy group $SU(3)$. This implies that $X$ has a holomorphic
volume form $vol_X \in H^{3,0}(X)$, unique up to scaling. Such a pair 
$(X,vol_X)$ is called a {\bf gauged}  Calabi Yau 3-fold. 
The {\bf (Kuranishi) moduli space} $S$ of $X$ is smooth and 
can be identified with a neighborhood of zero in $H^{2,1}(X) \cong H^1(X,T)$, 
s.\ \cite{Bo}, \cite{Ti} and \cite{To}. Denote by ${\cal X} = {(X_s)}_{s\in S}
 \rightarrow S$, $X_0 = X$, the corresponding deformation of complex 
structure. The {\bf ``intersection'' form} 
\begin{equation} \label{intersectionEqu} \omega (\xi , \eta ) 
:= \int_X \xi \wedge \eta\, ,\quad  
\xi , \eta \in H^3(X,\mbox{\Z})\, ,
\end{equation}
defines an integral nondegenerate skew symmetric bilinear form on
$H^3(X,\mbox{\Z})$. The corresponding complex symplectic form on
$H^3(X,\mbox{\C})$ will be denoted by the same letter. Consider 
the holomorphic line bundle $H^{3,0}({\cal X}) \rightarrow S$ with fibre 
$H^{3,0}(X_s)$ at $s\in S$. Denote by $H^{3,0}({\cal X}) - S$ the 
$\mbox{\C}^*$-bundle over $S$ which is obtained from the complex line 
bundle $H^{3,0}({\cal X})$ by removing the zero section $S \ni s \mapsto
0 \in H^{3,0}(X_s)$. We think of it as {\bf moduli space of gauged 
Calabi Yau 3-folds} $(X_s,vol_s)$,  $vol_s \in H^{3,0}(X_s) -\{ 0 \}$, 
$s\in S$. 

The holomorphic vector bundle $H^3({\cal X},\mbox{\C}) \rightarrow S$ has 
a canonical flat connection defined by the lattice bundle 
$H^3({\cal X},\mbox{\Z}) \subset H^3({\cal X},\mbox{\C})$, which is 
known as Gau{\ss}-Manin connection. Since the moduli space $S$ is local,
we can assume that $S$ is simply connected and that the bundle 
$H^3({\cal X},\mbox{\C}) \rightarrow S$ is trivial. In particular,
we have canonical identifications $H^3(X_s,\mbox{\Z}) \cong H^3(X,\mbox{\Z})$
and  $H^3(X_s,\mbox{\C}) \cong H^3(X,\mbox{\C})$. So we can define the 
{\bf period map} 
\[ {Per}: S \rightarrow P(H^3(X,\mbox{\C}))\, , \quad s\mapsto H^{3,0}(X_s)
\, .\] 
It follows from Kodaira and Spencer's deformation theory (s.\ e.g.\ \cite{M-K})
that
\[ d{Per}(T_sS) \; = \; d\pi (H^{3,0}(X_s) + H^{2,1}(X_s))\, ,\]
where $\pi :H^3(X,\mbox{\C}) \rightarrow P(H^3(X,\mbox{\C}))$ is the canonical
projection. This implies that the period map is an immersion and since 
$S$ is local we can assume that ${Per}: S \rightarrow Per(S) \subset 
P(H^3(X,\mbox{\C}))$ is an isomorphism. As a consequence, the cone 
$M_X = \cup_{s\in S}Per(s) -\{ 0\} \subset H^3(X,\mbox{\C})$ over 
$Per(S) = P(M_X)$ is canonically identified with the moduli space 
$H^{3,0}({\cal X}) - S$  of gauged Calabi Yau 3-folds. 
By the first Hodge-Riemann bilinear relations the tangent
space 
\[ T_uM_X \; = \;  H^{3,0}(X_s) + H^{2,1}(X_s)\, , \quad u\in Per(s)-\{ 0\}
\, ,\]
is a Lagrangean subspace of $H^3(X,\mbox{\C})$ with respect to the 
intersection form $\omega$. Remark that, using the standard real structure
$\tau$ on $H^3(X,\mbox{\C})$ with fix point set $H^3(X,\mbox{\R})$, we can
define the Hermitian form $\gamma$ of complex signature $(n+1,n+1)$, 
$n = h^{2,1}(X)$, on  $H^3(X,\mbox{\C})$ by equation (\ref{gammaEqu}). 
The first and second Hodge-Riemann bilinear relations imply that $M_X$ is a 
pseudo K\"ahler submanifold $M_X \subset (H^3(X,\mbox{\C}),\gamma )$ 
of complex signature $(1,n)$. More precisely, we have the relations: 
$\gamma (u,u) > 0$, $\gamma (v,v) < 0$ and $\gamma (u,v) = 0$ for all 
$u\in H^{3,0}(X_s) -\{ 0\}$ and $v\in H^{2,1}(X_s) -\{ 0\}$. This motivates 
the following definition.
\begin{dof} \label{formalDef} Given fundamental algebraic data 
$(V, \omega , \tau )$ as on 
p.~\pageref{fadP} and $\gamma$ defined in (\ref{gammaEqu}), a Lagrangean 
pseudo K\"ahler submanifold $M \subset (V, \omega ,\gamma )$ is called a 
{\bf formal moduli space (of gauged Calabi Yau 3-folds)} if the following
conditions are satisfied: 
\begin{itemize}
\item[(i)] $M$ is a cone, i.e.\ $\lambda M = M$ for all 
$\lambda \in \mbox{\C}-\{ 0\}$. 
\item[(ii)] $\gamma (u,u) > 0$ for all $u\in M$. 
\item[(iii)] $\gamma (v,v) < 0$ for all $0 \neq v \in T_uM$ such that 
$\gamma (u,v) = 0$.  
\end{itemize}   
\end{dof}  
 
To a Calabi Yau 3-fold $X$ we have associated 
the following algebraic data: $V = H^3(X,\mbox{\C})$, $\omega$ the  
intersection form (\ref{intersectionEqu}), $\tau$ the real structure with 
fix point set $V^{\tau} = H^3(X,\mbox{\R})$  and $\gamma = 
\sqrt{-1} \omega (\cdot , \tau \cdot )$. With this understood the next 
proposition gives a summary of the preceding discussion. 

\begin{prop} \label{formalProp} The cone $M_X \subset (V, \omega , \gamma )$ 
over the image 
$Per(S) = P(M_X)$ of the period map is a formal moduli space of gauged Calabi 
Yau 3-folds in the sense of Def.~\ref{formalDef}.  
\end{prop}      

\begin{prop} \label{exformalProp} The class of Lagrangean pseudo K\"ahler
cones defined in the previous section on p.\ \pageref{i-iiiDef} and classified
in \cite{dW-VP} and \cite{C2} consists of formal moduli spaces of gauged 
Calabi Yau 3-folds. 
\end{prop} 

Remark that any statement which is true for formal moduli spaces $M$ is true
when $M = M_X (\cong H^{3,0}({\cal X}) - S)$ is the (actual) moduli
space of gauged Calabi Yau 3-folds associated to the Kuranishi moduli space 
$S$ of a Calabi Yau 3-fold $X$. 

\noindent
{\bf Remark 5:} For any formal moduli space $M \subset (V,\omega , \gamma )$ 
there is a canonical K\"ahler metric $-g^{P(M)}$ on $P(M)$ known as 
{\bf special  K\"ahler metric}, cf.\ p. \pageref{sKP} and \cite{C2}, which 
can be defined by 
\begin{equation} \label{sKEqu} g^{P(M)}_{\pi u} (d\pi v, d\pi v) = 
\frac{\gamma (v,v)}{\gamma (u,u)} - 
\left|\frac{\gamma (u,v)}{\gamma (u,u)}\right|^2\, , \end{equation}  
for $u\in M$, $v\in T_uM$, where $\pi :M \rightarrow P(M)$ is the 
canonical projection. In the case of actual moduli spaces of gauged
Calabi Yau 3-folds $M = M_X$ the metric $g^{P(M)}$ is known as 
{\bf Weil-Petersson metric}.  The formal moduli spaces $M$ of 
Prop.~\ref{exformalProp} provide all the known examples of 
{\bf homogeneous} special K\"ahler manifolds, i.e.\ 
with transitive isometry group. It is noteworthy that ``most'' of these  
examples are not Hermitian symmetric, s.\ \cite{dW-VP} and \cite{C2}. 
We may ask the following natural question:\\
{\it Which of the homogeneous special K\"ahler manifolds can be realized 
as moduli spaces of Calabi Yau 3-folds (equipped with the Weil-Petersson 
metric)?}

To round up our presentation, we place the concept of formal moduli space
in the context of infinitesimal variations of Hodge structure, s.\ \cite{G}
and \cite{B-G}. Given fundamental algebraic data $(V, \omega, \tau )$ as 
on p.~\pageref{fadP}, a {\bf polarized Hodge structure of weight 3} on 
$V$, $\dim_{\mbox{\C}}V = 2n + 2$, with {\bf Hodge numbers} $h^{3,0} = 1$
and $h^{2,1} = n$ is given by a decomposition into complex subspaces  
\begin{equation}\label{HodgeEqu} V \; = \; H^{3,0} \oplus H^{2,1} 
\oplus H^{1,2} \oplus H^{0,3}\, ,
\end{equation} 
such that $H^{p,q} = \tau H^{q,p}$, $h^{p,q} = \dim_{\mbox{\C}}H^{p,q}$ and 
satisfying the first and second Hodge-Riemann bilinear relations. Moreover, 
one assumes that
a lattice $\Gamma \subset V^{\tau}$ is given and that $\omega$ restricts to an
{\it integral} nondegenerate skew symmetric bilinear form on $\Gamma$. 
Let us denote by $\cal D$ the classifying space for such Hodge structures
on $V$.  Remark that an element of $\cal D$, i.e.\ a Hodge structure as above,
is determined by the pair $(H^{3,0},H^{2,1})$. 

To any formal moduli space $M \subset  (V, \omega ,\gamma )$, s.\ 
Def.~\ref{formalDef}, we associate a map $M \rightarrow {\cal D}$ by 
\[ u \; \mapsto \; (H^{3,0}(u),H^{2,1}(u))\, , \]
\[ H^{3,0}(u)\; := \; \mbox{\C} u \, , \quad H^{2,1}(u) \; := \; 
T_uM \cap u^{\perp} \, ,\] 
where the orthogonal complement $\perp$ is to be taken with respect to the 
Hermitian form $\gamma$. Obviously this map factorizes to a map
$\varphi_M : P(M) \rightarrow {\cal D}$. If we denote by $p:  {\cal D}  
\rightarrow P(V)$ the projection $(H^{3,0},H^{2,1}) \mapsto H^{3,0}$
then $p \circ \varphi_M : P(M) \rightarrow P(V)$ is the trivial inclusion. 
For any Hodge decomposition (\ref{HodgeEqu}) one has also the corresponding
{\bf Hodge filtration} 
\[ F^3 \subset F^2 \subset F^1 \subset F^0 = V \, ,
\quad F^p \; = \; \oplus_{p\le k \le 3} H^{k,3-k}\, .\]
Consider now a variation of Hodge structure $u \mapsto (H^{p,q}(u))_{p,q}$, 
$u\in U$,  and 
denote by $(F^p(u))_p$ the corresponding Hodge filtrations. 
If the variation 
of Hodge structure arises from a local deformation of complex structure 
$(X_u)_{u\in U}$, i.e.\ $H^{p,q}(u) = H^{p,q}(X_u)$, then it must satisfy
Griffiths' infinitesimal period relations \cite{G}
\begin{equation} \label{iprEqu} \partial F^p(u) \; \subset \; F^{p-1}(u)\, .
\end{equation}  
This means that holomorphic partial derivatives of a (local) holomorphic 
section of the vector bundle $(F^p(u))_{u\in U}$ are sections of 
$(F^{p-1}(u))_{u\in U}$. The next theorem follows from the work of 
Bryant and Griffiths \cite{B-G}. 

\begin{thm}\label{moduliThm} Let  $M \subset (V, \omega ,\gamma )$ 
be a formal moduli space,
$\dim_{\mbox{\C}}V = 2n + 2$ and $\cal D$ the classifying space for 
Hodge structures as above. Then the map  $\varphi_M : P(M) 
\rightarrow {\cal D}$ is a solution to the differential system on $\cal D$  
defined by the infinitesimal period relations (\ref{iprEqu}). 
Conversely, any solution $\varphi : U \rightarrow {\cal D}$, $U$ a complex 
n-fold, to this differential system for which $p\circ \varphi : U 
\rightarrow P(V)$ is an immersion is locally of the form $\varphi_M$. 
\end{thm} 

Let $X$ be a Calabi Yau 3-fold, $S$ its Kuranishi moduli space and 
$M_X = \cup_{s\in S} Per(s) -\{ 0\} \subset H^3(X,\mbox{\C})$ the cone over
the image of the period map $Per: S \rightarrow P( H^3(X,\mbox{\C}))$. 
Recall that $M_X$ is the moduli space of gauged Calabi Yau 3-folds 
associated to $X$. The {\bf intermediate Jacobian} of $X_s$, $s\in S$, is 
the complex torus 
\[ {\cal J}(X_s) \; = \; \frac{H^3(X,\mbox{\C})}{H^{3,0}(X_s) + H^{2,1}(X_s) + 
H^3(X,\mbox{\Z})}\, .\] 
The {\bf bundle of intermediate Jacobians} ${\cal J} \rightarrow M_X$ over 
$M_X$ is the holomorphic torus bundle whose fibre at $u\in Per(s) -\{ 0\} \subset M_X$
is ${\cal J}_u = {\cal J}(X_s)$. 

\begin{thm} \label{JacThm} For any Lagrangean splitting 
$H^3(X,\mbox{\R}) = L_0 \oplus 
L_0'$ such that 
$M_X \subset H^3(X,\mbox{\C})$ is in general position, there is a pseudo hyper
K\"ahler structure of complex signature $(2,2n)$, $n = h^{2,1}(X)$, on    
the bundle of intermediate Jacobians ${\cal J} \rightarrow M_X$. 
\end{thm}

\noindent
{\bf Proof:} By Prop.\ \ref{formalProp} the moduli space of gauged Calabi
Yau 3-folds $M_X$ is a formal moduli space. In particular, s.\ 
Def.~\ref{formalDef}, it is a Lagrangean pseudo K\"ahler submanifold of 
$(H^3(X,\mbox{\C}),\omega , \gamma )$, where $\omega$ is the intersection form
(\ref{intersectionEqu}) and $\gamma = \sqrt{-1} \omega (\cdot , \tau \cdot )$
is defined with the help of the standard real structure $\tau$, i.e.\ 
$H^3(X,\mbox{\C})^{\tau} = H^3(X,\mbox{\R})$. 
By Thm.\ \ref{basicThm}, we can associate to these data together with the 
Lagrangean splitting of $H^3(X,\mbox{\R})$  a pseudo hyper K\"ahler structure
on $T^*M_X$.  Using the intersection form $\omega$ we can identify the 
cotangent bundle $T^*M_X \rightarrow M_X$ with the normal bundle ${\cal N}
\rightarrow M_X$ of $M_X \subset H^3(X,\mbox{\C})$. The normal bundle 
has fibre 
\[ {\cal N}_u \; = \; \frac{H^3(X,\mbox{\C})}{T_uM_X} \; = \; 
\frac{H^3(X,\mbox{\C})}{H^{3,0}(X_s) + H^{2,1}(X_s)}\] 
at $u\in Per(s) -\{ 0\} \subset M_X$ and ${\cal J} \rightarrow M_X$ is precisely 
the torus bundle ${\cal N}/ \Gamma \rightarrow M$ considered in Thm.\ 
\ref{torusThm} with $M = M_X$, $V = H^3(X,\mbox{\C})$ and 
$\Gamma = H^3(X,\mbox{\Z})$. 
Now the theorem is an immediate consequence of Thm.\ \ref{torusThm}. $\Box$ 

\noindent
{\bf Remark 6:} As checked in the proof of Prop.\ \ref{groupProp},
the standard complex symplectic structure $\Omega$ on $T^*M_X\cong {\cal N}$ is
invariant under the action of the lattice $\Gamma = H^3(X,\mbox{\Z}) 
\subset H^3(X,\mbox{\R}) = V^{\tau}$. So it factorizes to a complex 
symplectic structure on ${\cal J} = {\cal N}/\Gamma \cong T^*M/\Gamma$ 
independent of any choice of Lagrangean splitting $H^3(X,\mbox{\R})
= L_0\oplus L_0'$. This is the parallel complex symplectic structure
associated to the pseudo hyper K\"ahler structure of Thm.\ \ref{JacThm} 
and it coincides with the complex symplectic structure which was recently 
constructed by Donagi and Markman \cite{D-M}. 

\noindent
{\bf Remark 7:} Given a Lagrangean subspace $L_1 \subset (V,\omega ,\tau )$ 
such that $L_1 \cap \tau L_1 = \{ 0 \}$, there exists a Lagrangean
subspace $L$ such that $L= \tau L$ and $L\cap L_1 = \{ 0\}$. This shows that
locally one can always find a Lagrangean splitting of $H^3(X,\mbox{\R})$ such 
that $M_X$ is in general position. However, such a choice is not unique. 
In fact, the Lagrangean splittings of $H^3(X,\mbox{\R})\cong \mbox{\R}^{2n+2}$
are parametrized
by the coset space $Sp(n+1,\mbox{\R}) / GL(n+1,\mbox{\R})$. To reduce this 
arbitrariness 
to only a countable number of allowed choices, we can proceed as follows.   
Any isomorphism $H^3(X,\mbox{\Z}) \cong \mbox{\Z}^{n+1}
\oplus (\mbox{\Z}^{n+1})^*$ mapping the intersection form to the standard
nondegenerate integral skew symmetric bilinear form on $\mbox{\Z}^{n+1}
\oplus (\mbox{\Z}^{n+1})^*$ induces a Lagrangean splitting 
of $H^3(X,\mbox{\R})$. In other words, we allow only Lagrangean splittings of 
$H^3(X,\mbox{\R})$ which are induced by the choice of a symplectic 
basis $(\xi^i,\eta_i)_{i=0,\ldots , n}$ for the integral cohomology, i.e.\ 
$\xi^i,\eta_i \in H^3(X,\mbox{\Z})$, $\omega (\xi^i,\xi^j) = \omega 
(\eta_i,\eta_j) = 0$ and $\omega (\xi^i,\eta_j) = \delta^i_j$. 
The Lagrangean splittings 
of this type are parametrized by $Sp(n+1,\mbox{\Z})/GL(n+1,\mbox{\Z})$. 

\noindent
{\bf Remark 8:} The pseudo hyper K\"ahler metrics $G$ 
on the cotangent bundle $T^*M$ of a formal moduli space of gauged Calabi Yau
3-folds $M\subset V$, s.\ Def.\ \ref{formalDef} and Thm.\ \ref{basicThm}, 
are not complete. In fact, the open line segment $l$ joining a point $u\in M$
to the origin $0 \in V-M$ is a geodesic arc of finite length in the 
totally geodesic zero section $M \subset T^*M$. We may consider the natural
blow up $\sigma : \tilde{M} \rightarrow M$ of the cone $M$ at the origin. 
$\tilde{M} = M \cup P(M)$ is identified with the universal bundle   
of the projectivized cone $P(M)$. It is easy to see that 
$\sigma^*(G|M) = \sigma^* g$ extends smoothly to the divisor $P(M)\subset 
\tilde{M}$. However, this extension gives only a degenerate metric on 
$\tilde{M}\subset T^*\tilde{M}$, unless $M$ is a complex line. 
(Of course the same
remarks apply in the case of actual moduli spaces of gauged Calabi Yau
3-folds $M = M_X$ as considered in Thm.\ \ref{JacThm}.)

\end{document}